\documentstyle{book}
\include {mak}
\def\DJ{\leavevmode\setbox0=\hbox{D}\kern0pt\rlap
{\kern.04em\raise.188\ht0\hbox{-}}D}
\def\dj{\leavevmode
 \setbox0=\hbox{d}\kern0pt\rlap{\kern.215em\raise.46\ht0\hbox{-}}d}

\def\txt#1{{\textstyle{#1}}}
\def\hf{{\textstyle{1\over2}}}
\def\b{\beta}
\def\d{{\,\rm d}}
\def\e{\varepsilon}

\def\G{\Gamma}
\def\k{\kappa}
\def\s{\sigma}
\def\t{\theta}
\def\={\;=\;}

\def\zt{\zeta(\hf+it)}

\def\D{\Delta}

\def\R{\Re{\rm e}\,} \def\I{\Im{\rm m}\,}
\def\z{\zeta}

 \def\t{\theta}
\def\hf{{\textstyle{1\over2}}}
\def\txt#1{{\textstyle{#1}}}

\font\teneufm=eufm10
\font\seveneufm=eufm7
\font\fiveeufm=eufm5
\newfam\eufmfam
\textfont\eufmfam=\teneufm
\scriptfont\eufmfam=\seveneufm
\scriptscriptfont\eufmfam=\fiveeufm
\def\mathfrak#1{{\fam\eufmfam\relax#1}}

\font\tenmsb=msbm10
\font\sevenmsb=msbm7
\font\fivemsb=msbm5
\newfam\msbfam
\textfont\msbfam=\tenmsb
\scriptfont\msbfam=\sevenmsb
\scriptscriptfont\msbfam=\fivemsb
\def\Bbb#1{{\fam\msbfam #1}}

\def \NN {\Bbb N}

\def \RR {\Bbb R}
\def \ZZ {\Bbb Z}

\def\D{\Delta}

\def\b{\beta} \def\e{\varepsilon}
 \def\d{\,{\rm d}}

\begin{document}


\oddsidemargin 16.5mm
\evensidemargin 16.5mm

\thispagestyle{plain}

\noindent {\small\sc Univ. Beograd. Publ. Elektrotehn. Fak.}

\noindent {\scriptsize Ser. Mat. xx (xxxx), x--x.}

\vspace{5cc}
\begin{center}

{\Large\bf  ON SOME MEAN SQUARE ESTIMATES IN THE RANKIN-SELBERG PROBLEM
\rule{0mm}{6mm}\renewcommand{\thefootnote}{}%
}

\vspace{1cc}
{\large\it Aleksandar Ivi\'c}

\vspace{1cc}
\parbox{24cc}{{\scriptsize

An overview of the classical Rankin-Selberg problem involving the
asymptotic formula for sums of coefficients of holomorphic cusp
forms is given. We also study the function $\D(x;\xi)\;(0\le\xi\le1)$,
the error term in the
Rankin-Selberg problem weighted by $\xi$-th power of the logarithm.
Mean square estimates for $\D(x;\xi)$ are proved.

}}
\end{center}

\vspace{1.5cc}

\begin{center}
{\bf 1. THE RANKIN-SELBERG PROBLEM}
\end{center}

The classical Rankin-Selberg problem consists
of the estimation of the error term function
$$
\D(x) \;:=\; \sum_{n\le x}c_n - Cx,\leqno(1.1)
$$
where the notation is as follows.
Let $\varphi(z)$ be a holomorphic cusp form of weight $\kappa$ with
respect to the full modular group $SL(2,\ZZ)$, and denote by
$a(n)$ the $n$-th Fourier coefficient of $\varphi(z)$ (see e.g.,
R.A.  Rankin  [15] for a comprehensive account).
We suppose that $\varphi(z)$ is a normalized eigenfunction for the Hecke
operators $T(n)$, that is,
$  a(1)=1  $ and $  T(n)\varphi=a(n)\varphi $ for every $n \in \NN$.
In (1.1) $C>0$ is a suitable constant (see e.g., [9]
for its explicit expression), and $c_n$ is the convolution function
defined by
$$
c_{n}=n^{1-\kappa}\sum_{m^2 \mid n}m^{2(\kappa-1)}
\left|a\Bigl({n\over m^2}\Bigr)\right|^2.
$$
The classical Rankin-Selberg  bound of 1939 is
$$
\D(x) = O(x^{3/5}),\leqno(1.2)
$$
hitherto unimproved. In their works, done independently,
R.A. Rankin [16] derives (1.2) from a general result of
E. Landau [11], while A. Selberg [17] states the result with no proof.
Although the exponent 3/5 in (1.2) represents one of the longest
standing records in analytic number theory, recently there have been
some developments in some other aspects of the Rankin-Selberg
problem. In this paper we shall present an overview of some of these
new results. In addition, we shall consider the weighted sum (the so-called
Riesz logarithmic means of order $\xi$), namely
$$
{1\over\G(\xi+1)}\sum_{n\le x} c_n\log^\xi\left({x\over n}\right)
:= Cx + \D(x;\xi)\qquad(\xi\ge0),\leqno(1.3)
$$
where $C$  is as in (1.1), so that $\D(x) \equiv \D(x;0)$. The effect of
introducing weights such as the logarithmic weight in (1.3) is that the
ensuing error term (in our case this is $\D(x;\xi)$) can be estimated
better than the original error term (i.e., in our case $\D(x;0)$). This
was shown by Matsumoto, Tanigawa and the author in [9], where it was proved
that
$$
\D(x;\xi) \;\ll_\e\; x^{(3-2\xi)/5+\e}\qquad(0\le \xi \le 3/2).\leqno(1.4)
$$
Here and later $\e$ denotes arbitrarily
small constants, not necessarily the same ones at each occurrence, while
$a \ll_\e b$ means that the constant implied by the $\ll$-symbol depends on $\e$.
When $\xi =0$ we recover (1.2) from (1.4), only with the extra `$\e$' factor present.
In this work we shall pursue the investigations concerning $\D(x;\xi)$,
and deal with mean square bounds for this function.

\bigskip
\begin{center}
{\bf 2. THE FUNCTIONAL EQUATIONS}
\end{center}
In view of (1.1) and (1.2) it follows that the generating Dirichlet series
$$
Z(s) \;:= \; \sum_{n=1}^\infty c_nn^{-s}\qquad(s = \s+it)\leqno(2.1)
$$
converges absolutely for $\s>1$. The arithmetic function $c_n$
is multiplicative and satisfies $c_n \ll_\e n^\e$.
Moreover, it is well known (see e.g., R.A.
Rankin [14], [15]) that $Z(s)$ satisfies for all $s$ the functional equation
$$
\G(s+\k-1)\G(s)Z(s) = (2\pi)^{4s-2}\G(\k-s)\G(1-s)Z(1-s),\leqno(2.2)
$$
which provides then the analytic continuation of $Z(s)$.
In modern terminology  $Z(s)$ belongs to the Selberg
class $\cal S$ of $L$-functions of degree four (see A. Selberg   [18]
and the survey paper of Kaczorowski--Perelli [10]). An important feature,
proved by G. Shimura [19] (see also A. Sankaranarayanan [16]) is
$$
Z(s)  = \z(s)\sum_{n=1}^\infty b_n n^{-s}
= \z(s)B(s),\leqno(2.3)
$$
where $B(s)$ is holomorphic for $\s>0$, $b_n \ll_\e n^\e$ (in fact
$\sum_{n\le x}b_n^2 \le x\log^Ax$ holds, too). It also satisfies the functional
equation
\begin{eqnarray*}
B(s)\D_1(s) &=& B(1-s)\D_1(1-s),\\
\D_1(s) &=&
\pi^{-3s/2}\G(\hf(s+\k-1))\G(\hf(s+\k))\G(\hf(s+\k+1)),
\end{eqnarray*}
and actually $B(s) \in\cal S$ with degree three. The decomposition (2.3)
(the so-called `Shimura lift') allows one to use, at least to some extent,
results from the theory of $\z(s)$ in connection with $Z(s)$, and hence
to derive results on $\D(x)$.

\medskip
\begin{center}
{\bf 3. THE COMPLEX INTEGRATION APPROACH}
\end{center}
A natural approach to the estimation of $\D(x)$, used by the author in [8],
is to apply the classical complex integration technique. We shall
briefly present this approach now.
On using  Perron's inversion formula (see e.g., the Appendix
of [3]), the residue theorem and the convexity bound $Z(s) \ll_\e
|t|^{2-2\s+\e}\;(0 \le \s \le 1,\; |t| \ge 1)$, it follows that
$$
\D(x) = {1\over2\pi i}\int_{{1\over2}-iT}^{{1\over2}+iT}{Z(s)\over s}x^s\d s
+ O_\e\left(x^\e\left(x^{1/2} + {x\over T}\right)\right)\quad(1\ll T\ll x).\leqno(3.1)
$$
If we suppose that
$$
\int_X^{2X}|B(
\hf+it)|^2\d t \ll_\e X^{\t+\e}\qquad(\t\ge 1),\leqno(3.2)
$$
and use the elementary fact (see [3] for the results on the moments of $|\zt|$) that
$$
\int_X^{2X}|\zt|^2\d t \ll X\log X,\leqno(3.3)
$$
then from (2.3),(3.2),(3.3) and the Cauchy-Schwarz inequality for integrals
we obtain
$$
\int_X^{2X}|Z(\hf+it)|\d t  \ll_\e X^{(1+\t)/2+\e}.
$$
Therefore (3.1) gives
$$
\D(x) \ll_\e x^\e(x^{1/2}T^{\t/2-1/2} + xT^{-1}) \ll_\e x^{{\t\over
\t+1}+\e}\leqno(3.4)
$$
with $T = x^{1/(\t+1)}$. This was formulated in [8] as

\medskip
THEOREM A. {\it If $\t$ is given by } (3.2), {\it then}
$$
\D(x) \;\ll_\e \;x^{{\t\over
\t+1}+\e}.\leqno(3.5)
$$
\medskip

To obtain a value for $\t$,
note that  $B(s)$ belongs to the Selberg
class of degree three,
hence $B(\hf+it)$ in (3.2) can be written as a sum of two Dirichlet polynomials
(e.g., by the reflection principle discussed in [3, Chapter 4]),
each of length $\ll X^{3/2}$. Thus by
the mean value theorem for Dirichlet polynomials (op. cit.)
we have $\t \le 3/2$ in (3.2). Hence (3.5)  gives (with unimportant $\e$)
the Rankin-Selberg bound $\D(x) \ll_\e x^{3/5+\e}$.
Clearly improvement will come from better values of $\t$.
Note that the best possible
value of $\t$ in (3.2) is $\t = 1$, which follows
from general results on Dirichlet series
(see e.g., [3, Chapter 9]). It gives $1/2+\e$ as the
exponent in the Rankin-Selberg
problem, which is the limit of the method (the conjectural
exponent $3/8+\e$, which is best  possible,  is out of reach;
see the author's work [4]). To
attain this improvement one faces
essentially the same problem as in proving the sixth moment for $|\zt|$, namely
$$
\int_0^T|\zt|^6\d t \;\ll_\e \; T^{1+\e},
$$
only
this problem is even more difficult, because the
arithmetic properties of the coefficients $b_n$ are
even less known than the properties of
the divisor coefficients
$$
d_3(n) \= \sum_{abc=n;a,b,c\in\NN}1,
$$
generated by $\z^3(s)$.
If we knew the analogue of the strongest sixth
moment bound
$$
\int_0^T|\zt|^6\d t \;\ll\;T^{5/4}\log^CT\qquad(C>0),
$$
namely the bound (3.2) with $\t = 5/4$, then (3.1) would yield
$\D(x) \ll_\e x^{5/9+\e}$, improving substantially (1.2).

\medskip
The essential difficulty in this problem may be seen indirectly by comparing it
with the estimation of $\D_4(x)$, the error term in the asymptotic formula
for the summatory function of $d_4(n) = \sum_{abcd=n;a,b,c,d\in\NN}\,1$. The
generating function in this case is $\z^4(s)$. The problem analogous to the
estimation of $\D(x)$ is to estimate $\D_4(x)$, given the product representation
$$
\sum_{n=1}^\infty d_4(n)n^{-s} = \z(s)G(s)  = \z(s)\sum_{n=1}^\infty g(n)n^{-s}
\quad(\s>1)
\leqno(3.6)
$$
with $g(n) \ll_\e n^\e$ and $G(s)$ of degree three in the Selberg class
(with a pole of order three at $s=1$).
By the complex integration method one gets $\D_4(x) \ll_\e x^{1/2+\e}$ (here
$`\e$' may be replaced by a log-factor) using
the classical elementary bound $\int_0^T|\zt|^4\d t \;\ll\;T\log^4T$.
Curiously, this bound for $\D_4(x)$ has never been improved;
exponential sum techniques seem
to give a poor result here. However, if one knows only (3.6), then the situation
is quite analogous to the Rankin--Selberg problem, and nothing better than the
exponent 3/5 seems obtainable. The bound $\D(x) \ll_\e x^{1/2+\e}$ follows also
directly from (3.1) if the Lindel\"of hypothesis for $Z(s)$ (that
$Z(\hf+it) \ll_\e |t|^\e$) is assumed.

\begin{center}
{\bf 4. MEAN SQUARE OF THE RANKIN--SELBERG ZETA--FUNCTION}
\end{center}
\medskip
Let, for a given $\s\in\RR$,
$$
\mu(\s) = \limsup_{t\to\infty} {\log|\z(\s+it)|\over\log t}\leqno(4.1)
$$
denote the Lindel\"of function (the famous, hitherto unproved,
Lindel\"of conjecture for $\z(s)$ is that $\mu(\s) = 0$ for
$\s\ge\hf$, or equivalently that $\zt \ll_\e |t|^\e$). In [8] the
author proved the following

\bigskip
THEOREM B. {\it If $\b = 2/(5-\mu(\hf))$, then for fixed $\s$ satisfying
$\hf < \s \le 1$ we have}
$$
\int_1^{T}|Z(\s+it)|^2\d t = T\sum_{n=1}^\infty c_n^2n^{-2\s}
+ O_\e(T^{(2-2\s)/(1-\beta)+\e}).\leqno(4.2)
$$

\medskip
This result is the sharpest one yet when $\s$ is close to 1. For
$\s$ close to $\hf$ one cannot obtain an asymptotic formula, but only
the upper bound (this is [7, eq. (9.27)])
$$
\int_T^{2T}|Z(\s+it)|^2\d t \ll_\e\; T^{2\mu(1/2)(1-\s)+\e}(T+ T^{3(1-\s)})
\quad(\hf\le \s \le 1).\leqno(4.3)
$$
The upper bound in (4.3) follows easily from (2.3) and the fact that, as
already mentioned, $B(s)\in \cal S$ with degree three, so that $B(\hf+it)$
can be approximated by Dirichlet polynomials of length $\ll t^{3/2}$,
and the mean value theorem for Dirichlet polynomials yields
$$
\int_T^{2T}|B(\s+it)|^2\d t \ll_\e\; T^\e(T+ T^{3(1-\s)})
\qquad(\hf\le \s \le 1).
$$
Note that with the sharpest known result (see M.N. Huxley [2])
$\mu(1/2) \le 32/205$ we obtain $\b = 410/961 = 0.426638917\ldots\,$.
The limit is the value $\b = 2/5$ if the Lindel\"of hypothesis
(that $\mu(\hf) =0$) is true. Thus (4.2) provides a true asymptotic
formula for
$$
\s \;>\; {1+\b\over 2} \;=\; {1371\over1922} = 0.7133194\ldots\,.
$$
The proof of (4.2), given in [8], is based on the general method of the
author's paper [6], which contains a historic discussion on the formulas for the
left-hand side of (4.2) (see also K. Matsumoto [12]).

\medskip
We are able to improve (4.2) in the case when $\s = 1$. The result is
contained in

\medskip
THEOREM 1. {\it We have}
$$
\int_1^{T}|Z(1+it)|^2\d t = T\sum_{n=1}^\infty c_n^2n^{-2}
+ O_\e((\log T)^{2+\e}).\leqno(4.4)
$$

\medskip
{\bf Proof}. For $\s = \R s> 1$ and $X \ge 2$ we have
$$
\leqno(4.5)
$$
\begin{eqnarray*}
Z(s) &=& \sum_{n\le X}c_nn^{-s} + \int_X^\infty x^{-s}\d
\Bigl(\sum_{n\le x}c_n\Bigr)\\
&=&\sum_{n\le X}c_nn^{-s} + {CX^{1-s}\over s-1} - \D(x)X^{-s} -s\int_X^\infty
\D(x)x^{-s-1}\d x.
\end{eqnarray*}

By using (1.2) it is seen that the last integral converges absolutely for
$\s = \R s > 3/5$, so that (4.5) provides the analytic continuation of $Z(s)$
to this region. Taking $s =1+it, 1\le t\le T, X=T^{10}$, it follows that
$$
\int_1^{T}|Z(1+it)|^2\d t = \int_1^{T}
\Bigl\{\Bigl|\sum_{n\le X}c_nn^{-1-it}\Bigr|^2  -2C\I \Bigl(\sum_{n\le X}{c_n\over nt}
\Bigl({X\over n}\Bigr)^{it}\Bigr)\Bigr\}\d t + O(1).\leqno(4.6)
$$
By the mean value theorem for Dirichlet polynomials we have
\begin{eqnarray*}
\int_1^{T}\Bigl|\sum_{n\le X}c_nn^{-1-it}\Bigr|^2\d t
&=& T\sum_{n\le X}c_n^2 + O\Bigl(\sum_{n\le X}c_n^2n^{-1}\Bigr)\\
&=& T\sum_{n=1}^\infty c_n^2 + O_\e\Bigl((\log T)^{2+\e}\Bigr),
\end{eqnarray*}
where we used the bound (see K. Matsumoto [12])
$$
\sum_{n\le x}c_n^2 \;\ll_\e\; x(\log x)^{1+\e}\leqno(4.7)
$$
and partial summation. Finally we have
$$
\sum_{n\le X}{c_n\over n}\int_1^{T}{1\over t}
\left({X\over n}\right)^{it}\d t \ll \log\log T.\leqno(4.8)
 $$
 To see that (4.8) holds, note first that
for $X - X/\log T \le n \le X$ the integral over $t$ is trivially estimated
as $\ll \log T$, and the total contribution of such $n$ is
$$
\ll \log T\sum_{X - X/\log T \le n \le X}{c_n\over n}\d x \;\ll\; 1
$$
on using (1.1)--(1.2). For the remaining $n$ we note that the integral over $t$
equals
$$
{\left({X\over n}\right)^{it}\over it\log(X/n)}\Biggl|_1^T +
{1\over i\log(X/n)}\int_1^T\left({X\over n}\right)^{it}{\d t\over t^2}.
$$
The contribution of those $n$ is, using (1.1)--(1.2) again and making
the change of variable $X/u = v$,
\begin{eqnarray*}
&\ll&\sum_{1\le n\le X- X/\log T}{c_n\over n\log (X/n)} = \int_{1-0}
^{X-X/\log T}{1\over u\log (X/u)}\d(Cu + \D(u))\\
&=&\int_{1}^{X-X/\log T}{1\over u\log (X/u)}\left(C
+ {\D(u)\over u} + {\D(u)\over u\log(X/u)}\right)\d u + O(1)\\
&\ll&\int_{1}^{X-X/\log T}{\d u\over u\log (X/u)} + 1
= \int_{(1-1/\log T)^{-1}}^X{\d v\over v\log v} +1\\
&=& \log\log X - \log\log (1-1/\log T)^{-1} + 1 \ll \log\log T,
\end{eqnarray*}
and (4.8) follows.

 One can improve the error term in (4.4) to
$O(\log^2T)$, which is the limit of the method. I am very grateful
to Prof. Alberto Perelli, who has kindly indicated this to
me. The argument is very briefly as follows.
 Note that the coefficients $c_n^2$ are essentially the tensor product
of the $c_n$'s, and the $c_n$ are essentially the tensor product of
the $a(n)$'s; ``essentially" means in this case that the corresponding
$L$-functions differ at most by a ``fudge factor", i.e., a Dirichlet
series converging absolutely for $\sigma > 1/2$ and non-vanishing at
$s=1$. In terms of $L$-functions, the tensor product of the $a(n)$
(the coefficients of the tensor square $L$-function)
corresponds to the product of $\z(s)$ and the $L$-function of $Sym^2$
(Shimura's lift). Moreover, Gelbart--Jacquet [1] have shown that
$Sym^2$ is a cuspidal automorphic representation, so one can apply to
the above product the general Rankin-Selberg theory to obtain ``good
properties" of the corresponding $L$-function. Since $Sym^2$ is
irreducible, the $L$-function corresponding to $c_n^2$ has a double
pole at $s=1$ and a functional equation of Riemann type. It follows
that the sum in (4.7) is asymptotic to $Dx\log x$ for some $D>0$, and
the assertion follows by following the preceding argument.

\smallskip
In concluding this section, let it be mentioned that,
using (4.5), it easily follows that
$Z(1+it) \ll \log|t|\;(t\ge2)$.

\begin{center}
{\bf 5. MEAN SQUARE OF $\D(x;\xi)$}
\end{center}
\medskip
In this section we shall consider mean square estimates for $\D(x;\xi)$,
defined by (1.3). Although we could consider the range $\xi > 1$ as well,
for technical reasons we shall restrict ourselves to the range $0\le \xi \le1$,
which is  the condition that will be assumed henceforth to hold. Let
$$
\b_\xi \;:=\; \inf\,\Bigl\{\; \b \ge 0\;:\;
\int_1^X\D^2(x;\xi)\d x \ll X^{1+2\b}\;\Bigr\}.\leqno(5.1)
$$
The definition of $\b_\xi$ is the natural analogue of the classical
constants in mean square estimates for the generalized Dirichlet
divisor problem (see [3, Chapter 13]). Our first result in this direction is

\medskip
THEOREM 2. {\it We have}
$$
{3-2\xi\over8}\;\le\;\b_\xi\;\le\; \max\left({1-\xi\over2},\;{3-2\xi\over8}\,\right)
\qquad(0\le\xi\le 1).\leqno(5.2)
$$
\medskip
{\bf Proof}. First of all, note that (5.2) implies that $\b_\xi = (3-2\xi)/8$
for $\hf \le \xi\le 1$, so that in this interval the precise value of $\b_\xi$
is determined. The main tool in our investigations is the explicit Vorono{\"\i}
type formula for $\D(x;\xi)$. This is
$$
\D(x;\xi) \= V_\xi(x,N) + R_\xi(x,N),\leqno(5.3)
$$
where, for $N\gg1$,
$$
\leqno(5.4)
$$
\begin{eqnarray*}
V_\xi(x,N) &=& (2\pi)^{-1-\xi}x^{(3-2\xi)/8}\sum_{n\le N}c_nn^{-(5+2\xi)/8}
\cos\Bigl(8\pi(xn)^{1/4}+ \hf(\hf-\xi)\pi\Bigr),\\
R_\xi(x,N)  &\ll_\e& (xN)^\e\left(1 + x^{(3-\xi)/4}N^{-(1+\xi)/4}
+ (xN)^{(1-\xi)/4} + x^{(1-2\xi)/8}\right).
\end{eqnarray*}
This follows from the work of U. Vorhauer [20] (for $\xi=0$ this is
also proved in [9]), specialized to the case when
$$
A = {1\over(2\pi)^2}, B = (2\pi)^4, M = L = 2, b_1 = b_2 = d_1 = d_2 = 1,
\b_1 = \k - \hf, b_2 = \hf,
$$
$$
 \delta_1 = \k -{\txt{3\over2}}, \delta_2 = -\hf, \gamma = 1,
 p = B, q = 4, \lambda = 2, \Lambda = -1, C = (2\pi)^{-5/2}.
 $$
In (5.3)--(5.4) we take $N = x$, so that $R_\xi(x,N) \ll_\e x^{(1-\xi)/2+\e}$.
Since ${1-\xi\over2} \le {3-2\xi\over8}$ for $\xi \ge \hf$, the lower bound
in (5.2) follows by the method of [4]. For the upper bound we use $c_n\ll_\e n^\e$
and note that (${\rm e}(z) = \exp(2\pi iz))$
\begin{eqnarray*}
&&\int_X^{2X}\Bigl|\sum_{K<k\le2K}c_kk^{-(5+2\xi)/8}
{\rm e}(4(xk)^{1/4})\Bigr|^2 \d x\\&
\ll & X + \sum_{k_1\ne k_2}c_{k_1}c_{k_2}(k_1k_2)^{-(5+2\xi)/8}
\int_X^{2X}{\rm e}(4x^{1/4}(k_1^{1/4}-k_2^{1/4}))\d x\\
&\ll_\e& X + X^{3/4+\e}K^{-(5+2\xi)/4}\sum_{k_1\ne k_2}
{\Bigl|k_1^{1/4}-k_2^{1/4}\Bigr|}^{-1}\\
&\ll_\e& X + X^{3/4+\e}K^{(1-\xi)/2},
\end{eqnarray*}
where we used the first derivative test (cf. [3, Lemma 2.1]). Since $K \ll X$
and
$$
\int_X^{2X}\D^2(x;\xi)\d x \ll \int_X^{2X}|V_\xi(x,N)|^2\d x
+ \int_X^{2X}|R(x,N)|^2\d x,
$$
it follows that
$$
\int_X^{2X}\D^2(x;\xi)\d x \ll_\e X^{(7-2\xi)/4+\e} + X^{2-\xi+\e},
$$
which clearly proves the assertion.

\medskip
Our last result is a bound for $\b_\xi$, which improves on (5.2) when
$\xi$ is small. This is

\medskip
THEOREM 3. {\it We have}
$$
\b_\xi \;\le\; {2-2\xi\over5-2\mu(\hf)}\quad\qquad\Biggl(0\le\xi\le
{\txt{1\over6}}(1+2\mu(\hf))\Biggr).\leqno(5.5)
$$

\medskip
{\bf Proof}. We start from
$$
\D(x;\xi) = {1\over2\pi i}\lim_{T\to\infty}
\int_{c-iT}^{c+iT} Z(s){x^s\over s^{\xi+1}}\d s,\leqno(5.6)
$$
where $0 < c = c(\xi) < 1$ is a suitable constant (see K. Matsumoto [13]
for a detailed derivation of formulas analogous to (5.6)).
By the Mellin inversion formula we have (see e.g., the Appendix of [3])
$$
Z(s)s^{-\xi-1} = \int_0^\infty \D(1/x;\xi)x^{s-1}\d x \qquad(\R s = c).
$$
Hence by Parseval's formula for Mellin tranforms (op. cit.) we obtain, for
$\b_\xi < \s < 1$,
$$
\leqno(5.7)
$$
\begin{eqnarray*}
&&{1\over2\pi}\int_{-\infty}^\infty {|Z(\s+it)|^2\over
|\s+it|^{2\xi+2}}\d t = \int_0^\infty \D^2(1/x;\xi)x^{2\s-1}\d x\\
&=& \int_0^\infty \D^2(x;\xi)x^{-2\s-1}\d x \gg X^{-2\s-1}\int_X^{2X}\D^2(x;\xi)
\d x.
\end{eqnarray*}
Therefore if the first integral converges for $\s = \s_0+\e$, then (5.7) gives
$$
\int_X^{2X}\D^2(x;\xi) \;\d x\; \ll X^{2\s+1},
$$
namely $\b_\xi \le \s_0$. The functional equation (2.2) and Stirling's
formula in the form
$$
|\G(s)| = \sqrt{2\pi}|t|^{\s-1/2}{\rm e}^{-\pi|t|/2}(1+O(|t|^{-1}))
\qquad(|t| \ge t_0 > 0)
$$
imply that
$$
Z(s) \= {\cal X}(s)Z(1-s),\quad {\cal X}(\s+it) \asymp |t|^{2-4\s}
\quad(s = \s+it,\, 0\le \s \le 1, |t|\ge2).
\leqno(5.8)
$$
Thus it follows on using (4.3) that
\begin{eqnarray*}
\int_T^{2T}|Z(\s+it)|^2\d t &\ll& T^{4-8\s}\int_T^{2T}|Z(1-\s+it)|^2\d t\\
&\ll_\e& \,T^{4-8\s+2\mu({1\over2})\s + \max(1,3\s)+\e}.
\end{eqnarray*}
But we have
$$
4-8\s+2\mu(\hf)\s + \max(1,3\s) = 4 - 5\s + 2\mu(\hf)\s < 2\xi+2
$$
for
$$
\s > \s_0 = {2-2\xi\over5-2\mu(\hf)},\leqno(5.9)
$$
provided that $\s_0 \ge 1/3$, which occurs if $0 \le \xi \le{1\over6}(1+2\mu(\hf))$.
Thus the first integral in (5.7) converges if (5.9) holds, and Theorem 3 is proved.
Note that this result is a generalization of Theorem 7 in [8], which says that
$\b_0 \le (2-2\xi)/(5-2\mu(\hf))$.

\medskip
In the case when $\b_\xi = (3-2\xi)/8$ we could actually derive an asymptotic formula
for the integral of the mean square of $\D(x;\xi)$, much in the same way that
this was done in [9] for the square of
$$
\D_1(x) \;:=\; \int_0^x\D(u)\d u,
$$
where it was shown that
$$
\int_1^X\D_1^2(x)\d x = DX^{13/4} + O_\e(X^{3+\e})\leqno(5.10)
$$
with explicit $D>0$ (in [12] the error term was improved to $O_\e(X^3(\log X)^{3+\e})$).
In the case of $\D(x;1)$ the formula (5.10) may be used directly,
since
$$
{1\over x}\D_1(x) = {1\over x}\int_0^x\D(u)\d u = \D(x;1) + O_\e(x^\e).\leqno(5.11)
$$
To see that (5.11) holds, note that with $c=1-\e$ we have
\begin{eqnarray*}
\D(x;1) &=& {1\over2\pi i}\int_{c-i\infty}^{c+i\infty} Z(s){x^s\over s^2}\d s\\
&=& {1\over2\pi i}\int_{c-i\infty}^{c+i\infty} Z(s){x^s\over s(s+1)}\d s +
{1\over2\pi i}\int_{c-i\infty}^{c+i\infty} Z(s){x^s\over s^2(s+1)}\d s\\
&=& {1\over x}\int_0^x\D(u)\d u +
{1\over2\pi i}\int_{\e-i\infty}^{\e+i\infty} Z(s){x^s\over s^2(s+1)}\d s\\
&=& {1\over x}\D_1(x) +  O_\e(x^\e),
\end{eqnarray*}
on applying (5.8) to the last integral above.

\vspace{2cc}

\begin{center}{\small\bf REFERENCES}
\end{center}

\vspace{0.8cc}
\newcounter{ref}
\begin{list}{\small \arabic{ref}.}{\usecounter{ref} \leftmargin 4mm
\itemsep -1mm}



\item{} {\small {\sc S. Gelbart and H. Jacquet}, {\it A relation
between automorphic forms on
 GL{\rm(2)} and GL{\rm(3)}}, Proc. Nat. Acad. Sci. U.S.A. {\bf 73} (1976), 3348--3350.}

\item{} {\small {\sc M.N. Huxley}, {\it
Exponential sums and the Riemann zeta-function V},
Proc. London Math. Soc. (3) {\bf90}(2005), 1-41.}

\item{} {\small {\sc A. Ivi\'c}, {\it The Riemann zeta-function}, John Wiley \&
Sons, New York, 1985 (2nd ed., Dover, Mineola, N.Y., 2003).}

\item {} {\small  {\sc A. Ivi\'c}, {\it Large values of certain number-theoretic
error terms}, Acta Arith. {\bf56}(1990), 135-159.}

\item {} {\small  {\sc A. Ivi\'c}, {\it On some conjectures and results
for the Riemann zeta-function and Hecke series}, Acta Arith. {\bf99}(2001), 155-145.}

\item{} {\small {\sc A. Ivi\'c}, {\it On mean values
of zeta-functions in the critical strip},
J. Th\'eorie des Nombres de Bordeaux {\bf15}(2003), 163-173.}

\item{} {\small {\sc A. Ivi\'c},
{\it Estimates of convolutions of certain number-theoretic error terms},
Inter. J. of Math. and Mathematical Sciences
2004:1, 1-23.}

\item{} {\small {\sc A. Ivi\'c}, {\it
Convolutions and mean square estimates of certain number-theoretic
 error terms}, subm. to Publs. Inst. Math. (Beograd). {\tt arXiv.math.NT/0512306}}

\item{} {\small {\sc A. Ivi\'c, K. Matsumoto and Y. Tanigawa}, {\it On Riesz mean
of the coefficients of the Rankin--Selberg series}, Math. Proc.
Camb. Phil. Soc. {\bf127}(1999), 117-131.}

\item{} {\small {\sc A. Kaczorowski and A. Perelli}, {\it The Selberg class: a
survey}, in ``Number Theory in Progress, Proc. Conf. in honour
of A. Schinzel (K. Gy\"ory et al. eds)", de Gruyter,
Berlin, 1999, pp. 953-992.}

\item{} {\small {\sc E.~Landau}, {\it \"{U}ber die Anzahl der
Gitterpunkte in gewissen Bereichen II}, Nachr. Ges.
Wiss. G\"ottingen 1915, 209-243.}

\item{} {\small  {\sc K. Matsumoto},
{\it The mean values and the universality of Rankin-Selberg $L$-functions}.
M. Jutila (ed.) et al., Number theory.
``Proc. Turku symposium on number theory in memory of K. Inkeri",
Turku, Finland, May 31-June 4, 1999.  de Gruyter, Berlin, 2001, pp. 201-221.}

\item{} {\small {\sc K. Matsumoto}, {\it Liftings and mean
value theorems for automorphic $L$-functions}, Proc. London
Math. Soc. (3) {\bf90}(2005), 297-320.}

\item{} {\small {\sc R.A. Rankin}, {\it Contributions to the theory of Ramanujan's
function $\tau(n)$ and similar arithmetical functions. II, The order
of Fourier coefficients of integral modular forms},  Math. Proc.
Cambridge  Phil. Soc. {\bf 35}(1939), 357-372.}

\item{} {\small {\sc R.A. Rankin}, {\it Modular Forms}, Ellis Horwood Ltd.,
Chichester, England, 1984.}

\item{} {\small {\sc A. Sankaranarayanan}, {\it Fundamental properties of symmetric
square $L$-functions I}, Illinois J. Math. {\bf46}(2002), 23-43.}

\item{} {\small {\sc A.~Selberg}, {\it  Bemerkungen \"{u}ber eine Dirichletsche Reihe,
   die mit der Theorie der Modulformen nahe verbunden ist},
 Arch. Math. Naturvid. {\bf 43}(1940), 47-50.}

\item{}  {\small {\sc A.~Selberg}, {\it Old and new conjectures and results about a class
of Dirichlet series}, in ``Proc. Amalfi Conf. Analytic Number Theory 1989
(E. Bombieri et al. eds.)",
University of Salerno, Salerno, 1992, pp. 367--385.}

\item{} {\small {\sc G. Shimura}, {\it On the holomorphy
of certain Dirichlet series},
Proc. London Math. Soc. {\bf31}(1975), 79-98.}

\item{} {\small {\sc U. Vorhauer},
{\it Three two-dimensional Weyl steps in the circle problem II.
The logarithmic Riesz mean for a class of arithmetic functions}, Acta Arith.
{\bf91}(1999), 57-73.}
\vskip1cm
\end{list}
\vspace{1cc}


{\small
\noindent
Katedra Matematike RGF-a, Universitet u Beogradu,  \DJ u\v sina 7,
11000 Beograd, Serbia.

\noindent E--mail: ivic@rgf.bg.ac.yu, aivic@matf.bg.ac.yu

\end{document}